\documentclass{article}
\usepackage{amsfonts}
\usepackage{amsmath}
\usepackage{epsfig}

\newtheorem{theorem}{Theorem}

\newtheorem{example}[theorem]{Example}

\newenvironment{proof}[1][Proof]{\noindent\textbf{#1.} }{\ \rule{0.5em}{0.5em}}

\begin{document}

\title{A central limit theorem for time-dependent dynamical systems}
\author{P\'{e}ter N\'{a}ndori, Domokos Sz\'{a}sz and Tam\'{a}s Varj\'{u}
\thanks{
The support of the Hungarian National Foundation for Scientific
Research grant No. K 71693 is gratefully asknowledged.
}}

\maketitle

\begin{abstract}
The work  \cite{O-S-Y}  established memory loss
  in the time-dependent (non-random) case of uniformly expanding maps of the
   interval. Here we find conditions under which we have convergence to the normal distribution of the appropriately scaled
    Birkhoff-like partial sums of appropriate test functions. A substantial part of the problem is to ensure that the variances of the partial sums tend to infinity (cf. the zero-cohomology condition in the autonomous case). In fact, the present paper is the first one where non-random, i. e. specific examples are also found, which are not small perturbations of a given map.
 Our approach uses martingale approximation technique in the form of \cite{S-V}.
\end{abstract}

     \section{Introduction}
     Time-dependent dynamical systems appear in various applications. Recently,
     \cite{O-S-Y} could establish exponential loss of memory for expanding maps and, moreover,  for one-dimensional piecewise
     expanding maps with slowly varying parameters. It also provided interesting motivations and examples for the problem. For us - beside their work - an additional incentive was the question of J. Lebowitz \cite{Le}: bound the correlation decay for a planar finite-horizon Lorentz process which is periodic apart form the $0$-th cell; in it, the Lorentz particle encounters a particular scatterer of the $0$-th cell moderately  displaced at its each subsequent return to the $0$-th cell. (Slightly similar is the situation in the Chernov-Dolgopyat model of Brownian Brownian motion, where - between subsequent collisions of the light particle with the heavy one - the heavy particle slightly moves away, cf. \cite{ChD09}.)

The results of \cite{O-S-Y} say that - for sequences of uniformly uniformly expanding maps - distances of images of a pair of different initial measures converge to $0$ exponentially fast. In the same setup it is also natural to expect that probability laws of the Birkhoff-type partial sums of some given function - scaled, of course, by the square roots of their variances - are approximately Gaussian. The main theorem of our paper provides a positive answer though our conditions are surprisingly more restrictive than those of \cite{O-S-Y}. Let us explain the difficulty and some related results.

In functional central limit theorems for functions of autonomous
 chaotic  deterministic systems the zero-cohomology condition is - in quite a generality -  known to be necessary and sufficient for the vanishing of the limiting variance (see
   \cite{Li} for
    instance). For time-dependent systems, however, such a condition is only known for almost all versions of random
     dynamical systems  (see  \cite{A-L-S}) and for other models the situation can be and definitely is completely different. In fact, for time-dependent systems, first \cite{B} had proved a Gaussian approximation theorem in quite a generality; he, however, assumed that the variances of the Birkhoff-type partial sums tend to $\infty$ sufficiently fast; the paper, however, did not provide any example when this condition would hold.
     The more recent work \cite{CR} proves under some reasonable conditions a dichotomy: either the variances are bounded or the Gaussian approximation
     holds; the article also provides an example for the latter in the case when the time dependent maps are smaller and smaller perturbations of a given map. But still there is no general method for ascertaining whether the variance is bounded or not.
     Finally we note that \cite{KK} has interesting results for higher order cohomologies but its  setup is different.
     
 The present work is, in fact, the first one where non-random, i. e. specific examples are also found, that are not small perturbations of a given map. The proof of our main theorem uses martingale approximation technique in the form introduced in \cite{S-V} for treating additive functions of inhomogeneous Markov chains. The  organization of our paper is simple: its section 2 contains our main theorem and provides examples when it is applicable. Section 3 is devoted to the proof of the theorem.

\section{Results}

Let $A$ be a  set of numbers and $(X, \mathcal{F}, \mu)$ a probability
space. For each $a \in A$
define $T_{a}: X \rightarrow X$. Suppose that $\mu$ is invariant
for all $T_{a}$'s. Now consider a
sequence of numbers from $A$, i.e.
$\underline{a}: \mathbb{N} \rightarrow A$. Our aim is to
prove some kind of central limit theorem for the sequence
\[ f \circ T_{a_1}, f \circ T_{a_2} \circ T_{a_1}, ...\]
with some nice function $f: X \rightarrow \mathbb{R}$. \\
As usual,
\[ \hat{T}_{a}g (x) = g(T_{a}x)\]
and $\hat{T}^{*}$ is the $L^2(\mu)$-adjoint of $\hat{T}$
(the so called Perron-Frobenius
operator). Further, introduce the notation
$$
\hat{T}_{[i..j]} = \left\{ \begin{array}{rl}
  \hat{T}_{a_i}\dots \hat{T}_{a_j} &\mbox{ if $i \leq j$} \\
   Id &\mbox{ otherwise}
          \end{array} \right.
$$
and for simplicity write $\hat{T}_{[j]} = \hat{T}_{[1..j]}$.\\
Similarly, define
$$
\hat{T}_{[i..j]}^{*} = \left\{ \begin{array}{rl}
  \hat{T}_{a_j}^{*} \dots \hat{T}_{a_i}^{*} &\mbox{ if $i \leq j$} \\
     Id &\mbox{ otherwise}
               \end{array} \right.
$$
and $\hat{T}_{[j]}^{*} = \hat{T}_{[1..j]}^{*}$. \\
Further,
let $\mathcal{F}_0 = \mathcal{F}$,
$\mathcal{F}_i =  (T_{a_1})^{-1} \dots (T_{a_i})^{-1} \mathcal{F}_0 $ and
assume that there is a Banach space $\mathcal{B}$ of functions on $X$ such that
$\| g \| := \| g \|_{\mathcal{B}} \geq \| g\|_{\infty}$ for all
$g \in \mathcal{B}$. \\
Finally, for the fixed function $f$, introduce the notation
\[u_k = \sum_{i=1}^k \hat{T}_{[i+1..k]}^{*} f.\]
With the above notation, our aim is to prove limit theorem
for $S_n(x) = \sum_{k=1}^{n}  \hat{T}_{[k]} f(x)$.

\begin{theorem}
\label{tetel1}
Assume that $f$, $\underline{a}$ and $T_{b}$, $b\in A$
satisfy the following assumptions.
\begin{enumerate}
\item $\int f d \mu =0$.
\item $T_{b}$ is onto but not invertible for all $b \in A$.
\item $f \in \mathcal{B}$ and
there exist $K<\infty$ and $\tau <1$ such that
for all $\underline{b}$ sequences and for all $k$, $ \| \hat{T}_{b_1}^{*} ...
\hat{T}_{b_k}^{*} f \| < K \tau^k \| f \|$.

\item  \emph{(accumulated transversality)}
Define $\chi_k $ as the $L^2$-angle between
$u_k$ and the subspace of
$(T_{a_{k+1}})^{-1} \mathcal{F}_0$-measurable functions. Then
\[ \sum_{k=1}^N \min_{j \in \{ k, k+1\}} (1- \cos^2(\chi_j))\]
converges to $\infty$ as $N \rightarrow \infty$.
\end{enumerate}
Then
\[ Var(S_n) \rightarrow \infty\]
and
\[ \frac{S_n(x)}{\sqrt{ Var(S_n) }}\]
converges weakly to the standard normal distribution, where $x$ is distributed
according to $\mu$.
\end{theorem}

Assumption
3 roughly tells that there is an eventual spectral gap of the operators
$\hat{T}_{a_j}^{*}$ which is quite a natural assumption. Assumption 4 guarantees
that there is no much cancellation in $S_n$, for instance $f$ cannot be in
the cohomology class of the zero function when $|A|=1$.\\
Before proving the statement let us examine a special case.

\begin{example}
\label{pelda1}
Define $(X, \mathcal{F}, \mu) = (S^{1}$, Borel, Leb),
$A=\{ 2,3,\dots\}$, $T_{a}(x) = ax (mod 1)$, $\mathcal{B} = C^1=C^1(S^1)$,
\[\| g \| := \sup_{x \in S^1} |g(x)| + \sup_{x \in S^1} |g'(x)|. \]
Fix a non constant function $f \in C^1$ satisfying
$\int f d x=0$. Then there exists some integer $L=L(f)$ such that with all
sequences $\underline{a}$ for which
\[ \# \{ k: \min\{ a_k, a_{k+1}, a_{k+2} \} >L\} = \infty\]
the assumptions of Theorem 1 are fulfilled.
\end{example}

\begin{proof}[Proof of Example \ref{pelda1}]
It is easy to see that for all $g \in C^1$ with zero mean,
and for all $\underline{b} : \mathbb{N} \rightarrow A$,
\[ \| \hat{T}_{b_1}^{*}g \| \leq 2 b_{1}^{-1}\| g\| \]
and similarly,
\begin{equation}
\label{pelda_spectral_gap}
\| \hat{T}_{b_1}^{*} \dots \hat{T}_{b_k}^{*} g \| \leq 2 \cdot 2^{-k} \| g\|.
\end{equation}
Hence Assumption 3 is fulfilled.\\
In order to check Assumption 4, select $x,y \in S^{1}$, $\varepsilon >0, \delta >0$
such that
\[ \min_{z \in [ x, x+ \varepsilon]} f(z) >
\delta + \max_{z \in [ y, y+ \varepsilon]} f(z). \]
This can be done since $f$ is not constant.
Now choose
\[ L > \max\{\frac{16 \|f\|}{\delta} ,\frac{2}{ \varepsilon} \}.\]
Whence
\[ \| \hat{T}_{L}^{*}f \| \leq \delta /8. \]
Thus if $a_k > L$, then
\[ \| \sum_{i=1}^{k-1} \hat{T}_{[i+1..k]}^{*} f \| < 3 \delta/8 \]
is true independently of the choice of $a_1, \dots a_{k-1}$.
This yields
\[ \min_{z \in [ x, x+ \varepsilon]} u_k(z) >
\delta/4 + \max_{z \in [ y, y+ \varepsilon]}  u_k(z). \]
Since $L>\frac{2}{ \varepsilon}$, for all $g$ which is
$( T^{L})^{-1} \mathcal{F}_0$ measurable, one can find $h: [0, \varepsilon/2)
\rightarrow \mathbb{R}$ and $\varepsilon_1 \leq \varepsilon / 2$ such that
$g(y+\varepsilon_1+z) = g(x+z) = h(z)$ for all
$z \in [0, \varepsilon /2)$. Hence,
\begin{eqnarray}
&& \| u_k -g \|_2^2 \nonumber \\
&\geq& \int_{x}^{x + \varepsilon/2}
\left(  u_k(z) -g(z) \right)^2 dz
+\int_{y+ \varepsilon_1}^{x + \varepsilon_1 + \varepsilon/2}
\left(  u_k(z) -g(z) \right)^2 dz \nonumber \\
&=& \int_{0}^{\varepsilon/2}
\left(  u_k(x+z) -h(z) \right)^2 dz
+ \int_{0}^{\varepsilon/2}
\left(  u_k(y+ \varepsilon_1 +z) -h(z) \right)^2 dz \nonumber \\
&\geq& \frac{1}{2}  \int_{0}^{\varepsilon/2}
\left(  u_k(x+z) - u_k(y+ \varepsilon_1 +z) \right)^2 dz
\geq \frac{\delta^2 \varepsilon}{64} \label{mar23_1}
\end{eqnarray}
Since\[  \| u_k  \|_2 < \| u_k  \|\]
is bounded, (\ref{mar23_1}) implies that $(1-\cos^2(\chi_k)) $
is uniformly bounded away from zero if $\min \{ a_k, a_{k+1} \} >L$.\\
Hence, Assumption 5 is fulfilled if there exist infinitely many
indices $k$ such that
\[\min \{ a_k, a_{k+1}, a_{k+2} \} >L.\]
\end{proof}

In Example \ref{pelda1}, expanding maps with large derivative
were needed in order to obtain the Gaussian approximation. Naturally arises
the question that what happens in the case when one
uses only finitely many dynamics, for instance, only
$T_2$ and $T_3$ of Example \ref{pelda1}. That is why
we discuss the following example.

\begin{example}
\label{pelda2}
Define $X, \mathcal{F}, \mu, A, T_b, \mathcal{B}$ as in
Example \ref{pelda1}. If $\underline{a}$ is a sequence
for which there is a $b \in A$ such that for all integer $K$,
one can find a $k$ for which
\[ a_k = a_{k+1} = ... = a_{k+K-1} = b,\]
and $f \in \mathcal{B}$, $\int f =0$ is any function for which
the equation
$f = \hat{T}_b u - u$ has no solution $u$, then the assumptions of
Theorem 1 are fulfilled.
\end{example}

\begin{proof}[Proof of Example \ref{pelda2}]
It is enough to verify Assumption 4. To do so,
for $K \in \mathbb{Z}_{+}$ pick $k$ such that
\begin{equation}
\label{pelda2_1}
a_{k-K} = a_{k-K+1} = ... = a_{k+2} = b.
\end{equation} Then
(\ref{pelda_spectral_gap}) implies that
\begin{equation}
\label{pelda2_2}
 \| u_j - \sum_{i=0}^{\infty} \left( \hat{T}_b^{*} \right)^{i} f\|
< C 2^{-K}
\end{equation}
holds for $j=k,k+1$ with some $C$ uniformly in $K$. Now, if
$g :=  \sum_{i=0}^{\infty} \left( \hat{T}_b^{*} \right)^{i} f$
is not $(T_b)^{-1} \mathcal{F}_0$-measurable, then necessarily
its $L^2$-angle with those functions is positive.
Since (\ref{pelda2_1}) and (\ref{pelda2_2}) hold for infinitely
many $k$'s, $\min \{\chi_k, \chi_{k+1} \}$ has a positive lower
bound infinitely many times, inferring Assumption 5.
On the other hand, if $g$ is
$(T_b)^{-1} \mathcal{F}_0$-measurable, then $g = \hat{T}_b
\hat{T}_b^{*}g$
and $g- \hat{T}_b^{*}g = f$ imply that for $u= \hat{T}_b^{*}g$,
$\hat{T}_b u - u =f$.
\end{proof}

Note, that in Example \ref{pelda2}, $Var(S_n)$ can be arbitrary small.
Indeed, pick a $C^1$ function $f$, for which $f = \hat{T}_3 u - u$ has
no solution $u$, but there is some $v$ such that $f = \hat{T}_2 v - v$.
Now, pick a sequence of integers $d_l, l \in \mathbb{N}$,
$d_l \rightarrow \infty$ fast enough, and define
$$
a_k= \left\{ \begin{array}{rl}
  3&\mbox{ if $\exists l: d_l \leq k < d_l+l$} \\
     2 &\mbox{ otherwise.}
               \end{array} \right.
$$
It is easy to see that (\ref{pelda_spectral_gap}) implies
$\mathbb{E}(|\hat{T}_{[i]}f \cdot \hat{T}_{[j]}f|) \leq 2^{|i-j|+1} \|f \|^2$
(formally it follows from (\ref{exp_lecsenges})), which in turn yields
that $Var(S_k)$ is bounded by some constant times $k$.
Now, with the notation $l_n:= \max\{l: d_l \leq n \}$, write
\begin{eqnarray*}
 Var(S_n) &\leq& 4 Var(S_{d_{l_n-1}+l_n}) + 4 Var(S_{d_{l_n}}-
S_{d_{l_n-1}+l_n})\\ 
&&+ 4 Var(S_{d_{l_n}+l_n} - S_{d_{l_n}} )
+ 4 Var(S_n - S_{d_{l_n}+l_n}).
\end{eqnarray*}
On the other hand, $f = \hat{T}_2 v - v$ implies that
$\hat{T}_2 f + ... + \hat{T}_2^m f$
is uniformly bounded in $m$. Thus the second and the last term
in the above sum are bounded. Whence $Var(S_n)$ is smaller than
some constant times $d_{l_n-1}$.
Especially, if $d_l=2^{2^{2^l}}$, then
\[ \frac{Var(S_n)}{n^{\alpha}} \rightarrow 0\]
as $n \rightarrow 0$ for any $\alpha$ positive. Note that
in this case the conditions of \cite{B} for the
Gaussian approximation are not met.

\section{Proof of Theorem \ref{tetel1}}
This section is devoted to the proof of Theorem \ref{tetel1}.

As in \cite{Li}, \cite{S-V} and \cite{CR}, the proof is based on martingale approximation.
First, observe that
\[ \hat{T}_{[n]}^{*} \hat{T}_{[n]} = Id\]
and
\[  \hat{T}_{[n]} \hat{T}_{[n]}^{*} \]
is the orthogonal projection onto the $\mathcal{F}_n$ measurable functions (for the
proof of the latter, see \cite{Li}). Now we introduce our approximating
martingale, which is analogous to the one of \cite{S-V}:
\begin{eqnarray}
\label{Z_k_def}
Z_k = \sum_{i=1}^k \mathbb{E} \left[ \hat{T}_{[i]} f |
\mathcal{F}_k\right] = \sum_{i=1}^k \hat{T}_{[k]}\hat{T}_{[k]}^{*}
\hat{T}_{[i]}f
= \sum_{i=1}^k \hat{T}_{[k]} \hat{T}_{[i+1..k]}^{*} f =
\hat{T}_{[k]} u_k
\end{eqnarray}
Since
\begin{eqnarray}
\label{Zfkapcsolat}
\hat{T}_{[i]} f &=& Z_i - \mathbb{E} \left[ Z_{i-1} |
\mathcal{F}_i \right] \\
&=& \left( Z_i - \mathbb{E} \left[ Z_{i} | \mathcal{F}_{i+1} \right] \right)
+ \left( \mathbb{E} \left[ Z_{i} | \mathcal{F}_{i+1} \right] -
 \mathbb{E} \left[ Z_{i-1} | \mathcal{F}_{i} \right] \right),
\end{eqnarray}
one obtains
\begin{equation*}
S_n = \sum_{k=1}^{n-1}
\left( Z_k - \mathbb{E} \left[ Z_{k} | \mathcal{F}_{k+1} \right] \right)
+ Z_n.
\end{equation*}
Now,
\[ \xi_k^{(n)} = \frac{1}{\sqrt{Var (S_n)}}
\left( Z_k - \mathbb{E} \left[ Z_{k} | \mathcal{F}_{k+1} \right] \right),\]
is a reverse martingale difference for the $\sigma$-algebras
$\mathcal{F}_1, \dots \mathcal{F}_n$.
Thus, in particular
\begin{equation}
\label{varfelbontas}
Var (S_n) = \sum_{k=1}^{n-1} Var
\left( Z_k - \mathbb{E} \left[ Z_{k} | \mathcal{F}_{k+1} \right] \right)
+ Var(Z_n).
\end{equation}
Using our martingale approximation and the well known martingale CLT
(see \cite{S-V} for instance), it is enough to prove that the difference
between the martingale approximant and $S_n$ is negligible,
\begin{equation}
\label{CHT_1feltetel}
\max_{1 \leq i \leq n} \| \xi_i^{(n)} \|_{\infty} \rightarrow 0
\end{equation}
and
\begin{equation}
\label{CHT_2feltetel}
\| \sum_{i=1}^n \mathbb{E} \left[ \left(  \xi_i^{(n)}\right)^2
| \mathcal{F}_{i+1} \right] - 1 \|_2 \rightarrow 0.
\end{equation}
To prove (\ref{CHT_1feltetel}) and (\ref{CHT_2feltetel}),
we adopt the ideas of \cite{S-V}.
To verify (\ref{CHT_1feltetel}), observe that by Assumption 4,
\begin{eqnarray}
\label{Z_k_becsles}
\| Z_k \|_{\infty} &\leq& \sum_{j=1}^k \|  \hat{T}_{[k]}  \hat{T}_{[j+1..k]}^{*}
f \|_{\infty}
\leq \sum_{j=1}^k \|  \hat{T}_{[j+1..k]}^{*} f \|_{\infty} \nonumber \\
&\leq& \sum_{j=1}^k \|  \hat{T}_{[j+1..k]}^{*} f \|
\leq \sum_{j=1}^k K \tau^{k-j} \| f \| \leq C_f. \label{Z_k_becsles}
\end{eqnarray}
Thus
\begin{equation}
\label{Z_k_becsles_2}
 \| \mathbb{E} \left[ Z_{k} | \mathcal{F}_{k+1} \right] \|_{\infty}
\leq C_f.
\end{equation}
Now, we prove that the variance of $S_n$ converges to infinity:
\begin{equation}
\label{var_becsles}
Var(S_n) = \mu (S_n^2) \rightarrow \infty
\end{equation}
as $n \rightarrow \infty$.
Since (\ref{Z_k_becsles}) implies that $Var(Z_n)$ is bounded,
(\ref{varfelbontas}) can be written as
\begin{eqnarray*}
Var(S_n)&=& O(1) + \sum_{k=1}^{n-1} \mathbb{E}(Z_k^2) +
\mathbb{E} \left( \mathbb{E}[Z_k|\mathcal{F}_{k+1}]^2 \right) -
2 \mathbb{E} \left( Z_k \mathbb{E}[Z_k|\mathcal{F}_{k+1}] \right) \\
&=& O(1) + \sum_{k=1}^{n-1} \mathbb{E}(Z_k^2) -
\mathbb{E} \left( \mathbb{E}[Z_k|\mathcal{F}_{k+1}]^2 \right)\\
&=& O(1) +
\sum_{k=1}^{n-1} \| u_k \|_2^2 - \| u_k \|_2^2 \cos^2 \chi_k.
\end{eqnarray*}
Here, we used (\ref{Z_k_def}), and the fact that $\hat{T}_{[k]}$ is
$L^2(\mu)$-isometry.
Now, since
\[ Var(f) = Var(\hat{T}_{[i]} f) \leq 2 Var(Z_i) +
2 Var(\mathbb{E}[Z_{i-1}|\mathcal{F}_{i}])
\leq 2 \| u_i \|_2^2 + 2 \| u_{i-1} \|_2^2,\]
one obtains
\[ Var(S_n) \geq O(1) + \frac{1}{4} Var(f)  \sum_{k=1}^{n-1}
\min_{j \in \{ k, k+1\}} \left( 1- \cos^2 \chi_j\right),\]
which converges to infinity as $n \rightarrow \infty$ by
Assumption 4. Thus we have verified (\ref{var_becsles}). \\
Now, (\ref{Z_k_becsles}), (\ref{Z_k_becsles_2}) and
(\ref{var_becsles}) together imply (\ref{CHT_1feltetel})
and that
the difference between the martingale and $S_n$ is negligible.\\
To verify (\ref{CHT_2feltetel}), first observe that for $i > j$
\begin{eqnarray}
 \| \mathbb{E} \left[ \hat{T}_{[j]} f | \mathcal{F}_i \right]\|_{\infty}
&=&  \| \hat{T}_{[i]} \hat{T}_{[i]}^{*}  \hat{T}_{[j]} f \|_{\infty}
= \| \hat{T}_{[i]} \hat{T}_{[j+1..i]}^{*} f \|_{\infty}
= \|  \hat{T}_{[j+1..i]}^{*} f \|_{\infty} \nonumber \\
& \leq& K \tau^{i-j} \| f \|. \label{exp_lecsenges}
\end{eqnarray}
Then one can prove the assertion obtained from Lemma 4.4 in \cite{S-V}
by replacing $v_l^{(n)}$ with
\[ \mathbb{E} \left[ \left( \xi_{n-l}^{(n)} \right)^2
| \mathcal{F}_{n-l+1} \right]\]
the same way as it was done in \cite{S-V}, which yields (\ref{CHT_2feltetel}).

\section*{Acknowledgements}
The authors are highly indebted to Mikko Stenlund and Lai Sang Young for first explaining their result in October 2010 and second for a most valuable discussion in April 2011.


\begin{thebibliography}{99}

\bibitem{A-L-S}  Ayyer, A., Liverani, C., Stenlund, M.: Quenched CLT for
Random Toral Automorphisms, Discrete and Continuous Dynamical
Systems, 24 331-348. (2009)

\bibitem{B} V. I. Bakhtin, Random processes generated by a hyperbolic sequence of mappings. I, Russian
Acad. Sci. Izv. Math. 44 (1995), no. 2, 247-279,
Random processes generated by a hyperbolic sequence of mappings. II, Russian Acad. Sci.
Izv. Math. 44 (1995), no. 3, 617-627.

\bibitem{ChD09} Chernov, N., Dolgopyat. D.: Brownian Brownian Motion--1, Memoirs AMS. {\bf 198}, No. 927, pp 193.

\bibitem{CR} Conze, J. P., Raugi, A.: Limit theorems for sequential expanding dynamical  systems of [0,1],
Contemporary Mathematics, 430, 89-121 (2007)

\bibitem{KK} A. Katok, S. Katok. Higher cohomology for abelian groups of toral automorphisms,  Ergod. Th. \& Dynam. Sys.  I. {\bf 15}, 569-592, 1995;
II.  {\bf 25}, 1909-1917  (2005)

\bibitem{Le} Lebowitz, J. L. Oral communication, (2005)

\bibitem{Li} Liverani, C.: Central Limit Theorem for Deterministic Systems,
International Conference on Dynamical Systems, Montevideo 1995,
Pitman Research Notes in Mathematics Series, 362 (1996)

\bibitem{O-S-Y} Ott, W., Stenlund, M., Young, L.-S.: Memory Loss for
Time-Dependent Dynamical Systems. Math. Res. Lett. 16,
463-475. (2009)

\bibitem{S-V} Sethuraman, S., Varadhan, S. R. S.: A Martingale Proof of
Dobrushin's Theorem for Non-Homogeneous Markov Chains, Electronic J. Prob. 10.
(2005)

\end{thebibliography}
\end{document}